\documentclass[a4paper,11pt,twoside]{article}

\RequirePackage[a4paper,head=1cm]{geometry}
\usepackage[utf8x]{inputenc}
\usepackage{amsmath}
\usepackage{amssymb}
\usepackage{hyperref}
\usepackage{dsfont}

\newtheorem{Theoreme}{Theorem}
\newtheorem{Lemme}{Lemma}[section]

\newtheorem{Proposition}{Proposition}[section]
\newtheorem{Remarque}{\bf Remark}

\title{\bf Some functional inequalities on \\  polynomial volume growth Lie groups} 
\author{Diego Chamorro } 

\begin{document} 

\maketitle 
\begin{center}
\begin{minipage}[l]{140mm}
\begin{scriptsize}\abstract{We study in this article some Sobolev-type inequalities on polynomial volume growth Lie groups. We show in particular that improved Sobolev inequalities can be extended without the use of the Littlewood-Paley decomposition to this general framework.\\[3mm] 
\textbf{Keywords:} Sobolev inequalities, polynomial volume growth Lie groups.}
\end{scriptsize}
\end{minipage}
\end{center}

\section{Introduction}
Classical Sobolev inequalities provide us with a family of a priori estimates of the following form
\begin{equation}\label{ClassicalSobo}
\|f\|_{L^q}\leq C\|\nabla f\|_{L^p} \quad \mbox{where } q=np/(n-p).
\end{equation}
Initially stated over $\mathbb{R}^n$, they were succesively generalized to other settings such as manifolds or Lie groups, see for example  \cite{Badr} or \cite{Varopoulos} for such generalizations. 

Since the work of P. Gérard, Y. Meyer \& F. Oru in \cite{GMO}, we know that it is possible to improve the classical Sobolev inequalities (\ref{ClassicalSobo}) by introducing a well-suited Besov space, and it is worth knowing if these improved inequalities can be generalized to some special Lie groups. For example, in the case of the Heisenberg group -which is the simplest example of stratified Lie groups- this was done by H. Bahouri, P. Gérard \& C-J Xu in \cite{BGX} following essentially the same ideas of the original paper of P. Gérard, Y. Meyer and F. Oru; while, for general stratified Lie groups, the task was achieved in \cite{CHAME1} using some different techniques. In this special setting we obtained a family of Sobolev-type inequalities:  namely, for $G=(\mathbb{R}^n, \cdot, \delta)$ a stratified Lie group and for $f$ a function such that $ f\in \dot{W}^{s,p}(G)$ and $f\in \dot{B}^{-\beta, \infty}_{\infty}(G)$, we have:
\begin{equation}\label{ImprovedClassicalSobo}
\|f\|_{\dot{W}^{s,q}}\leq C \|f\|_{\dot{W}^{s_1,p}}^\theta\|f\|_{\dot{B}^{-\beta,\infty}_{\infty}}^{1-\theta}
\end{equation}
where the parameters $p,q$ and the indexes $\theta, \beta, s $ and $s_1$ are related in a specific manner. See section \ref{FSD} below for the definitions of these functional spaces.\\

This type of Lie groups is a generalization of $\mathbb{R}^n$ when modifying dilations; and in this setting, mathematical objects we are dealing with are constructed in such way in order to respect the homogeneity induced by these dilations. Therefore, many properties of these objects (operators, functional spaces) are very similar to the Euclidean case. See \cite{Folland2} and  \cite{Stein2} and the references given there for more details. 

If we want to study inequalities of type (\ref{ClassicalSobo}) and (\ref{ImprovedClassicalSobo}) in a more abstract framework, it is possible to consider Lie groups without a dilation structure and in this case we have several possibilities: a first example is given by \textit{nilpotent Lie groups} which are a generalization of stratified Lie groups (recall that every stratified Lie group is nilpotent) but these groups are not necessarily endowed with a dilation structure, see more details in \cite{Folland}. A second example is given by \textit{polynomial volume growth Lie groups}, where we have useful polynomial estimates for the Haar measure of a ball. Some other examples can be considered such as \textit{exponential growth Lie groups}, see the book \cite{Varopoulos} for definitions and some related results for the last case.

Classical Sobolev inequalities have been extensively studied in the three previous frameworks and a detailed account can be found in \cite{Varopoulos}.\\

In this article we will especially focus on \textit{polynomial volume growth Lie groups} and our main purpose is to study improved Sobolev inequalities of type (\ref{ImprovedClassicalSobo}) in this very particular setting.\\ Our results are as follows:

\begin{Theoreme}\label{smile00} 
Let $G$ be a polyomial volume growth Lie group.  
\begin{enumerate}
\item[$\bullet$] \textbf{\emph{[Strong inequalities $p>1$]}}
If $ f\in \dot{W}^{s_1,p}(G)$ and $f\in \dot{B}^{-\beta, \infty}_{\infty}(G)$, then we have
\begin{equation}\label{PGLR1}
\|f\|_{\dot{W}^{s,q}}\leq C \|f\|_{\dot{W}^{s_1,p}}^\theta\|f\|_{\dot{B}^{-\beta,\infty}_{\infty}}^{1-\theta}
\end{equation}
where $1<p<q<+\infty$, $\theta=p/q$, $s=\theta s_1 -(1-\theta)\beta$ and $-\beta<s<s_1$.

\item[$\bullet$] \textbf{\emph{[Strong inequalities $p=1$]}}
If $ \nabla f\in L^1(G)$ and $f\in \dot{B}^{-\beta, \infty}_{\infty}(G)$, then we have
\begin{equation}\label{smile0} 
\|f\|_ {L^{q}}\leq C \|\nabla f \|_ {L^{1}}^{\theta} \|f\|_{\dot{B}^{-\beta, \infty}_{\infty}}^{1-\theta} 
\end{equation} 
where $1<q<+\infty$, $\theta = 1/q$ and $\beta=\theta/(1-\theta)$.  
\item[$\bullet$] \textbf{\emph{[Weak inequalities $p=1$]}} 
If $ \nabla f\in L^1(G)$ and $f\in \dot{B}^{-\beta, \infty}_{\infty}(G)$, then we have
\begin{equation}\label{smile1} 
\|f\|_ {\dot{W}^{s, q}_{\infty}}\leq C \|\nabla f\|_{L^1}^{\theta}\|f\|_{\dot{B}^{-\beta,\infty}_{\infty}}^{1-\theta} 
\end{equation} 
where $1<q<+\infty$, $0<s<1/q<1$, $\theta=1/q$ and $\beta=\frac{1-sq}{q-1}$.
\end{enumerate}
\end{Theoreme} 
For a precise definition of these functional spaces, refer to section \ref{FSD}.\\

Let us make some remarks concerning the techniques used in the proof of theses inequalities. For the proof of inequality (\ref{PGLR1}), we will not make use of the Littlewood-Paley decompostion as it is done in \cite{GMO} or \cite{BGX} and this estimate will be obtained in a more straightforward way applying the sub-Laplacian's fractional powers properties together with some properties of the Besov spaces.

Concerning strong inequalities (\ref{smile0}) and weak inequalities (\ref{smile1}), the proof follows a very different path and we will see that these two estimates rely on the Modified Poincaré pseudo-inequality stated in theorem \ref{CHAME}. Observe in particular that it is the use of this special inequality that suggest us the definition of the weak Sobolev spaces $\dot{W}^{s, q}_{\infty}$ in the estimates (\ref{smile1}).

To finish preliminary remarks, let us stress here that the role of the polynomial growth geometry will be clearly identified in the estimates available for the heat kernel $H_t$ associated to the sub-Laplacian $\mathcal{J}$ and in the operator's properties build from the spectral decomposition of the sub-Laplacian $\mathcal{J}$. These estimates and properties will be decisive for the proof of theorem \ref{smile00}.\\

The plan of the article is the following: section \ref{PVGLG} is devoted to a short introduction of polynomial volume growth Lie groups, section \ref{LHK} gives some important estimations for the Heat kernel, section \ref{SDSL} presents some results concerning spectral theory, section \ref{FSD} gives the precise definition of functional spaces involved in the inequalities above, while section \ref{proofs} presents the proof of theorem \ref{smile00}. 

\section{Polynomial volume growth Lie groups}\label{PVGLG}
Let $G$ be a connected unimodular Lie group endowed with its Haar measure $dx$.  Denote by $\mathfrak{g}$ the Lie algebra of $G$ and consider a family ${\bf X}=\{X_1,...,X_k\}$ of left-invariant vector fields on $G$ satisfying the Hörmander condition, which means that the Lie algebra generated by the $X_j$ for $1\leq j\leq k$ is $\mathfrak{g}$.\\

In this setting we have at our disposal the Carnot-Carathéodory metric associated with ${\bf X}$ defined as follows: let $\ell:[0,1]\longrightarrow G$ be an absolutely continuous path. We say that $\ell$ is admissible if there exists measurable functions $\gamma_1,...,\gamma_k:[0,1]\longrightarrow\mathbb{C}$ such that, for almost every $t\in [0,1]$, we have $\ell'(t)=\sum_{j=1}^k\gamma_j(t)X_j(\ell(l))$. If $\ell$ is admissible, define the length of $\ell$ by $\|\ell\|=\int_{0}^1(\sum_{j=1}^k|\gamma_j(t)|^2)^{1/2}dt$. Then, for all $x,y\in G$, the distance between $x$ and $y$ is the infimum of the lengths of all admissible curves joining $x$ to $y$. We will denote $\|x\|$ the distance between the origin $e$ and $x$ and $\|y^{-1}\cdot x\|$ the distance between $x$ and $y$.\\

For $r>0$ and $x\in G$, denote by $B(x,r)$ the open ball with respect to the Carnot-Carathéodory metric centered in $x$ and of radius $r$, and by $V(r)$ the Haar measure of any ball of radius $r$. When $0<r<1$, there exists $d\in \mathbb{N}^\ast$, $c_l$ and $C_l>0$ such that, for all $0<r<1$ we have
$$c_l r^d\leq V(r)\leq C_l r^d.$$
The integer $d$ is the \textit{local} dimension of $(G, {\bf X})$. When $r\geq 1$, two situations may occur, independently of the choice of the family ${\bf X}$: either $G$ has polynomial volume growth and there exist $D\in \mathbb{N}^\ast$,  $c_\infty$ and $C_\infty>0$ such that, for all $r\geq 1$ we have
\begin{equation}\label{GLCP}
c_\infty r^D\leq V(r)\leq C_\infty r^D,
\end{equation}
or $G$ has exponential volume growth, which means that there exist  $c_e, C_e, \alpha, \beta>0$ such that, for all $r\geq 1$ we have
$$c_e e^{\alpha r}\leq V(r)\leq C_e  e^{\beta r}.$$
When $G$ has polynomial volume growth, the integer $D$ in (\ref{GLCP}) is called the dimension at infinity of $G$. Recall that nilpotent groups have polynomial volume growth and that a strict subclass of the nilpotent groups consists of stratified Lie groups. For more details see the book \cite{Varopoulos}.\\

We will assume from now on that $G$ is a connected unimodular polynomial volume Lie group with local dimension $d$ and dimension at infinity $D$.

\section{Sub-Laplacian and Heat kernel}\label{LHK}
Once we have fixed the family $\textbf{X}$, we define the gradient on $G$ by $\nabla = (X_1,...,X_k)$ and we consider a sub-Laplacian $\mathcal{J}$ on $G$ defined by $\mathcal{J}=-\sum_{j=1}^k X_j^2$, which is a positive self-adjoint, hypo-elliptic operator since ${\bf X}$ satisfies the Hörmander's condition. Its associated heat operator on $G\times]0, +\infty[$ is given by $\partial_{t}+\mathcal{J}$. We recall in the next theorem some well-known properties of the semi-group $H_t$ obtained from the sub-Laplacian $\mathcal{J}$. See the book \cite{Varopoulos} and the references given there for a proof.
\begin{Theoreme}\label{fol} There exists a unique family of continuous linear operators $(H_{t})_{t>0}$ defined on $L^{1}+L^{\infty}(G)$ with the semi-group property $H_{t+s}=H_{t}H_{s}$ for all $t, s>0$ and $H_{0}=Id$, such that: 
\begin{enumerate} 
\item[1)] the sub-Laplacian $\mathcal{J}$ is the infinitesimal generator of the semi-group  $H_{t}=e^{-t\mathcal{J}}$;  
\item[2)] $H_{t}$ is a contraction operator on $L^{p}(G)$ for $1\leq p\leq +\infty$ and for $t>0$; 
\item[3)] the semi-group $H_t$ admits a convolution kernel $H_{t}f=f\ast h_{t}$ where $h_{t}$ is the heat kernel.
\item[4)] $ \|H_{t}f-f \|_ {L^p}\to 0$ if $t\to 0$ for $f \in L^{p}(G)$ and $1\leq p < +\infty$; 
\item[5)] If $f\in L^{p}(G)$, $1\leq p\leq +\infty$, then the function $u(x, t)=H_{t}f(x)$ is a solution of the heat equation.
\end{enumerate} 
\end{Theoreme} 
We obtain in particular that $H_t$ is a symmetric diffusion semi-group as considered by Stein in \cite{Stein0} with infinitesimal generator $\mathcal{J}$. \\

We need to fix some terminology. To begin with note that associated to the family $\textbf{X}$ we also have a family of right-invariant vector fields $\{Y_1,...,Y_k\}$ with similar properties. Let $I=(j_{1},...,j_{\beta})\in \{1,...,k\}^{\beta}\; (\beta \in \mathbb{N})$ be a multi-index, we set $|I|=\beta$ and define $X^I$ and $Y^I$ by the formula $X^{I}=X_{j_{1}}\cdots X_{j_{\beta}}$ ($Y^{I}=Y_{j_{1}}\cdots Y_{j_{\beta}}$ resp.) with the convention $X^{I}=Id$ if $\beta=0$. The interaction of operators $X^I$ and $Y^I$ with convolutions is clarified by the following identities:
$$ X^I(f \ast g) = f \ast (X^I g), \quad   Y^I (f \ast g) = (Y^I f ) \ast g,   \quad  (X^I f ) \ast g = f  \ast (Y^I g).$$
In particular we have $(\nabla f)\ast g= f \ast(\tilde{\nabla} g)$ where $\tilde{\nabla}=(Y_1,...,Y_k)$.\\

We will say now that $\varphi\in \mathcal{C}^{\infty}(G)$ belongs to the Schwartz class $\mathcal{S}(G)$ if 
$$N_{\alpha,I}(\varphi)=\underset{x\in G}{ \sup}(1+\|x\|)^{\alpha}|X^{I}\varphi(x)|<+\infty. \qquad(\alpha\in \mathbb{N}, I\in \underset{\beta\in \mathbb{N}}{\cup} \{1,...,k\}^{\beta}).$$

\begin{Remarque}
\emph{To characterize the Schwartz class $\mathcal{S}(G)$ we can replace vector fields $X^I$ in the semi-norms $N_{k,I}$ above by right-invariant vector fields $Y^I$.}
\end{Remarque}
For a proof of these facts and for further details see \cite{Folland2}, \cite{Stein0}, \cite{Varopoulos} and the references given there. 

\begin{Theoreme}\label{TheoCalor}
Let $G$ be a polynomial volume growth Lie group, then for every $j\in \{1,...,k\}$, there exists $C>0$ such that
\begin{equation*}
\left|X_jh_{t}(x)\right|\leq C t^{-1/2}V(\sqrt{t})^{-1}e^{-\frac{\|x\|^2}{ct}}\qquad  \mbox{for all } x\in G, \; t>0.
\end{equation*}
\end{Theoreme}
This theorem implies the next proposition
\begin{Proposition}\label{poly5}
For every  $j\in \{1,...,k\}$ and for all $p\in [1,+\infty]$ there exists a constant $C>0$ such that:
\begin{equation}\label{uni8bis}
\|X_{j}h_{t}(\cdot)\|_{L^p}\leq C t^{-1/2}V(\sqrt{t})^{-1/p'},\quad t>0; 
\end{equation}
\end{Proposition}
For a proof of theorem \ref{TheoCalor} and proposition \ref{poly5}, see chapter VIII of the book \cite{Varopoulos}.

\section{Spectral decomposition for the sub-Laplacian}\label{SDSL}

The use in this article of spectral resolution for the sub-Laplacian consists roughly in expressing this operator by the formula $\mathcal{J}=\int_{0}^{+\infty}\lambda \;dE_{\lambda}$ and, by means of this characterization, build a family of new operators $m(\mathcal{J})$ associated to a Borel function $m$. This kind of operators have some nice properties as shown in the next propositions.

\begin{Proposition} If $G$ is polynomial growth Lie group and if $m$ is a bounded Borel function on $]0,+\infty[$ then the operator $m(\mathcal{J})$ defined by
\begin{equation}\label{little} 
m(\mathcal{J})=\int_{0}^{+\infty}m(\lambda) \;dE_{\lambda }, 
\end{equation} 
is bounded on $L^{2}(G)$  and admits a convolution kernel $M$ i.e.: $m(\mathcal{J})(f)=f\ast M \qquad (\forall f\in L^{2}(G))$. 
\end{Proposition} 

Following \cite{Hulanicki} and \cite{Furioli2} we can improve the conclusion of the above proposition. Let $k\in\mathbb{N}$ and $m$ be a function of class $\mathcal{C}^{k}(\mathbb{R}^{+})$, we write 
\begin{equation*} 
\|m \|_ {(k)}=\underset{\underset{\lambda>0}{1\leq r\leq k}}{\sup } (1+\lambda)^{k }|m^{(r)}(\lambda)|.
\end{equation*}
This formula gives us a necessary condition to obtain some properties of the operators defined by (\ref{little}):

\begin{Proposition}\label{toge} Let $G$ be polynomial volume growth Lie group with local dimension $d$. Let  $j\in \{1,...,k\}$ and $p\in [1, +\infty]$. There is a constant $C>0$ and an integer $k$ such that, for any function $m\in \mathcal{C}^{k}(\mathbb{R}^{+})$ with $ \|m\|_{(k)}<+\infty$, the kernel $M_{t}$ associated to the operator $m(t\mathcal{J})$ with $t>0$ satisfies 
 \begin{equation} \label{poly8}
\|X_jM_{t}(\cdot)\|_{L^p}\leq Ct^{-(\frac{d}{2p'}+\frac{1}{2})}\|m\|_{(k)}. 
\end{equation} 
where $\frac{1}{p}+\frac{1}{p'}=1$.
\end{Proposition} 
\textbf{\textit{Proof}}. Follow the same steps of the proof of proposition 3.2 in \cite{Furioli2} and use inequality (\ref{uni8bis}).
\begin{flushright}{$\blacksquare$}\end{flushright} 
\begin{Remarque}
\emph{Notice that, when $0<t\leq 1$, we can replace in (\ref{poly8}) $X_j$ by $X^I$ for some multi-index $I$.}
\end{Remarque}

\section{Functional spaces}\label{FSD} 

We give in this section the precise definition of the functional spaces involved in theorem \ref{smile00}. In a general way, given a norm $\|\cdot\|_{E}$, we will define the corresponding functional space $E(G)$ by $\{f\in \mathcal{S}'(G): \|f\|_{E}<+\infty\}$. For the Lebesgue spaces $L^p(G)$ with $1\leq p+\infty$, we will use the following characterization
\begin{equation*} 
\|f\|^{p}_{L^{p}}=\int_{0}^{+\infty}p\sigma^{p-1 } |\{x\in G:|f(x)|> \sigma\}|d\sigma,  
\end{equation*} 
and for the Lorentz spaces $L^{p,\infty}(G)$ we set $\|f\|_ { L^{p, \infty}}=\underset{\sigma>0}{\sup}\{\sigma \; |\{x\in G:|f(x)|> \sigma\}|^{1/p}\}$.\\

In order to define Sobolev spaces, we need to introduce the fractional powers $\mathcal{J}^s$ and $\mathcal{J}^{-s}$ with $s>0$:
\begin{eqnarray}
\mathcal{J}^s f(x)&=&\underset{\varepsilon \to 0}{\lim}\frac{1}{\Gamma(k-s)}\int_{\varepsilon}^{+\infty}t^{k-s-1}\mathcal{J}^k H_tf(x)dt\nonumber \\
\mathcal{J}^{-s} f(x)&=&\underset{\eta \to +\infty}{\lim}\frac{1}{\Gamma(s)}\int_{0}^{\eta}t^{s-1}H_tf(x)dt \label{fractionalJ}
\end{eqnarray}
for all $f\in \mathcal{C}^{\infty}(G)$ with $k$ the smallest integer greater than $s$. We consider then the Sobolev spaces by the norms
\begin{equation}\label{DefSobosp}
\|f\|_ {\dot{W}^{s, p}} =\|\mathcal{J}^{s/2}f\|_{L^{p}}
\end{equation}
when $1<p<+\infty$ and when $p=s=1$ we will note 
\begin{equation}\label{DefSobo1p}
\|f\|_ {\dot{W}^{1,1}} =\|\nabla f\|_{L^{1}}.
\end{equation}
We will also need to define weak Sobolev spaces $\dot{W}^{s,p}_\infty(G)$ used in (\ref{smile1}) and we write here
\begin{equation} \label{DefDebilSoboLe}
\|f\|_ {\dot{W}^{s, p}_\infty} =\|\mathcal{J}^{s/2}f\|_{L^{p,\infty}} \qquad (1<p<+\infty) 
\end{equation}
 
Finally, for Besov spaces of indices $(-\beta, \infty, \infty)$ which appear in all the inequalities (\ref{PGLR1})-(\ref{smile1}), we have:
\begin{equation}\label{BesovChaleurGhomo}
\|f\|_{\dot{B}^{-\beta,\infty}_{\infty}}=\underset{t>0}{\sup}\;\;t^{\beta/2 } \|H_{t}f\|_ {L^\infty} 
\end{equation} 
The choice of this \textit{thermic} definition for Besov spaces will be clarified in the next section. Observe that other equivalent characterizations do exist in the framework of polynomial volume growth Lie groups, see for example \cite{Furioli2} or \cite{Saka}, but they are not as useful in our computations as the thermic one.

\section{Improved Sobolev Inequalities on stratified groups: the proofs}\label{proofs} 

We will divide the proof of the theorem \ref{smile00} in two steps following the values of the parameter $p$ used in the Sobolev spaces which appear on the right hand side of inequalities (\ref{PGLR1})-(\ref{smile1}). This separation of the proof in the cases when $p>1$ and when $p=1$ is due to the definition of Sobolev spaces given by the formulas (\ref{DefSobosp}) and (\ref{DefSobo1p}) and is independent from the underlying geometry.
Thus, we first study the inequality (\ref{PGLR1}) and then we prove the strong inequality (\ref{smile0}) and the weak inequality (\ref{smile1}). 

\subsection{The general improved Sobolev inequalities ($p>1$)}
We start the proof observing that the operator  $\mathcal{J}^{s/2}$ carries out an isomorphism between the spaces $\dot{B}^{-\beta, \infty}_{\infty}(G)$ and $\dot{B}^{-\beta-s, \infty}_{\infty}(G)$. This fact follows from the thermic definition of Besov spaces (see \cite{Saka} for a proof and more details). We can rewrite the inequality (\ref{PGLR1}) in the following way 
\begin{equation*}
\| \mathcal{J}^{\frac{s-s_1}{2}}f\|_{L^q}\leq C \|f\|_{L^p}^\theta\|f\|_{\dot{B}^{-\beta-s_1,\infty}_{\infty}}^{1-\theta}
\end{equation*}
where $1<p<q<+\infty$, $\theta=p/q$, $s=\theta s_1 -(1-\theta)\beta$ and $-\beta<s<s_1$. Using the sub-Laplacian fractional powers characterization (\ref{fractionalJ}) we have the identity
\begin{equation}\label{Kashmor}
\mathcal{J}^{\frac{-\alpha}{2}}f(x)=\frac{1}{\Gamma(\frac{\alpha}{2})}\int_{0}^{+\infty}t^{\frac{\alpha}{2}-1}H_t f(x)dt=\frac{1}{\Gamma(\frac{\alpha}{2})}\left(\int_{0}^{T}t^{\frac{\alpha}{2}-1}H_t f(x)dt+\int_{T}^{+\infty}t^{\frac{\alpha}{2}-1}H_t f(x)dt\right)
\end{equation}
where $\alpha=s_1-s>0$ and where $T>0$ is a parameter that will be fixed in the sequel.\\

For studying each one of these integrals we will use the estimates
\begin{enumerate}
\item[$\bullet$] $|H_tf(x)|\leq |f(x)|$\\

\item[$\bullet$] $|H_tf(x)|\leq C t^{\frac{-\beta-s_1}{2}}\|f\|_{\dot{B}^{-\beta-s_1, \infty}_{\infty}}$ \qquad (by the thermic definition of Besov spaces)\\
\end{enumerate}
Then, applying these inequalities in (\ref{Kashmor}) we obtain
$$|\mathcal{J}^{\frac{-\alpha}{2}}f(x)|\leq \frac{c_1}{\Gamma(\frac{\alpha}{2})}T^{\frac{\alpha}{2}} |f(x)|+ \frac{c_2}{\Gamma(\frac{\alpha}{2})}T^{\frac{\alpha-\beta-s}{2}}\|f\|_{\dot{B}^{-\beta-s_1, \infty}_{\infty}}.$$
We fix now $$T=\left(\frac{\|f\|_{\dot{B}^{-\beta-s_1, \infty}_{\infty}}}{|f(x)|}\right)^{ \frac{2}{\beta+s_1}}$$
and we get 
$$|\mathcal{J}^{\frac{-\alpha}{2}}f(x)|\leq \frac{c_1}{\Gamma(\frac{\alpha}{2})}|f(x)|^{1-\frac{\alpha}{\beta+s_1}}+ \frac{c_2}{\Gamma(\frac{\alpha}{2})}|f(x)|^{1-\frac{\alpha}{\beta+s_1}}\|f\|^{\frac{\alpha}{\beta+s_1}}_{\dot{B}^{-\beta-s_1, \infty}_{\infty}}.$$
Since $\frac{\alpha}{\beta+s_1}=1-\theta$ and $\theta=p/q$ we have
$$|\mathcal{J}^{\frac{-\alpha}{2}}f(x)|\leq \frac{c}{\Gamma(\frac{\alpha}{2})}|f(x)|^{\theta}\|f\|^{1-\theta}_{\dot{B}^{-\beta-s_1, \infty}_{\infty}}.$$
We finally obtain
$$\|\mathcal{J}^{\frac{-\alpha}{2}}f\|_{L^q}\leq c\|f\|_{L^p}^{\theta}\|f\|^{1-\theta}_{\dot{B}^{-\beta-s_1, \infty}_{\infty}}$$
and we are done. 
\begin{flushright}{$\blacksquare$}\end{flushright} 

\subsection{Strong an weak inequalities ($p=1$)}\label{ISPGH} 
We treat now the inequalities (\ref{smile0}) and (\ref{smile1}). For this we will need the following result.

\begin{Theoreme}[Modified Poincaré pseudo-inequality]\label{CHAME} Let $f$ be a function such that $\nabla f\in L^1(G)$. We have the following estimate for $0 ≤s< 1$ and for $t > 0$:
\begin{equation}\label{MPPI}
\|\mathcal{J}^{s/2}f-H_{t}\mathcal{J}^{s/2}f\|_{L^{1}}\leq C \; t^{\frac{1-s}{2}}\|\nabla f\|_{ L^{1}}.  
\end{equation} 
\end{Theoreme}

Let us make some remarks. This theorem is crucial for proving strong and weak inequalities when $p=1$, mainly because this estimate is especially well suited for matching with the Besov space's thermic definition. Note also that, when $s=0$, we have some alternative proofs of 
(\ref{MPPI}) depending on the framework and its underlying geometry. See \cite{Ledoux} for details.  In the general case exposed in theorem \ref{CHAME}, the role of the geometry is given in the $L^p$-estimates available for the Heat kernel. \\

\emph{\textbf{Proof}}.
To begin the proof, we observe that the following identity occurs: 
\begin{equation*} 
(\mathcal{J}^{s/2}f-H_{t}\mathcal{J}^{s/2}f)(x)=\left(\int_{0}^{+\infty}m(t\lambda)dE_{\lambda}\right)t^{1-s/2}\mathcal{J} f(x), 
\end{equation*} 
where we noted $m(\lambda)=\lambda^{s/2-1}(1-e^{-\lambda }$) for $\lambda>0$, note that $m$ is a bounded function which tends to $0$ at infinity since $s/2-1<0$. We break up this function by writing:  
$$m(\lambda)=m_{0}(\lambda)+m_{1}(\lambda)=m(\lambda)\theta_{0}(\lambda)+m(\lambda)\theta_{1}(\lambda)$$ 
where we chose the auxiliary functions $\theta_{0}(\lambda), \theta_{1}(\lambda)\in \mathcal{C}^{\infty}(\mathbb{R}^{+})$ defined by:  
\begin{eqnarray*} 
\bullet\quad \theta_{0}(\lambda) = 1 \quad\mbox{on } \quad]0, 1/2 ] \quad\mbox{and}\quad 0 \quad \mbox{on}\quad]1, +\infty[, \\[5mm ] 
\bullet\quad \theta_{1}(\lambda) = 0 \quad\mbox{on} \quad]0, 1/2 ] \quad\mbox{and}\quad 1 \quad \mbox{on }\quad]1, +\infty[, 
\end{eqnarray*} 
so that $\theta_{0}(\lambda)+\theta_{1}(\lambda)\equiv 1$. Then, we obtain the formula: 
\begin{equation*} 
(\mathcal{J}^{s/2}f-H_{t}\mathcal{J}^{s/2}f)(x)=\left(\int_{0}^{+\infty}m_{0}(t\lambda)dE_{\lambda}\right)t^{1-s/2}\mathcal{J}f(x)+
\left(\int_{0}^{+\infty}m_{1}(t\lambda)dE_{\lambda}\right)t^{1-s/2}\mathcal{J}f(x).  
\end{equation*} 
If we note $M^{(i)}_{t}$ the kernel of the operator fixed by $\int_{0}^{+\infty}m_{i}(t\lambda)dE_{\lambda}$ for $i=0,1$, we have:  
\begin{equation*} 
(\mathcal{J}^{s/2}f-H_{t}\mathcal{J}^{s/2}f)(x)=t^{1-s/2}\mathcal{J}f\ast M^{(0)}_{t}(x)+t^{1-s/2}\mathcal{J}f\ast M^{(1)}_{t}(x).  
\end{equation*} 
We obtain the inequality
\begin{equation}\label{poids1} 
\int_{G}\left|\mathcal{J}^{s/2}f-H_{t}\mathcal{J}^{s/2}\right|dx \leq \int_{G}\left|t^{1-s/2}\mathcal{J}f\ast M^{(0)}_{t}(x)\right|dx +\int_{G}\left|t^{1-s/2}\mathcal{J}f\ast M^{(1)}_{t}(x)\right|dx.
\end{equation} 
We will now estimate the right side of the above inequality by the two following propositions:  
\begin{Proposition}\label{Sir} 
For the first integral in the right-hand side of (\ref{poids1}) we have the inequality:  
\begin{equation*} 
\int_{G}\left|t^{1-s/2}\mathcal{J}f\ast M^{(0)}_{t}(x)\right|dx\leq Ct^{\frac{1-s}{2}}\|\nabla f\|_{L^{1}} 
\end{equation*}
\end{Proposition} 
\emph{\textbf{Proof}}. The function $m_{0}$ is the restriction on $\mathbb{R}^{+}$ of a function belonging to the Schwartz class. This function satisfies the assumptions of the proposition \ref{toge} which we apply after having noticed the identity
\begin{equation*} 
I=\int_{G}\left|t^{1-s/2}\mathcal{J}f\ast M^{(0)}_{t}(x)\right|dx=\int_{G}\left|t^{1-s/2}\nabla f\ast\tilde{\nabla}M^{(0)}_{t}(x)\right|dx 
\end{equation*} 
where we noted $\tilde{\nabla}$ the gradient formed by the vectors fields $(Y_{j})_{1\leq j\leq k}$. We have then
\begin{equation*} 
I\leq \int_{G}\int_{G}t^{1-s/2}|\nabla f(y)||\tilde{\nabla} M^{(0)}_{t}(y^{-1}\cdot x)|dxdy\leq t^{1-s/2}\|\nabla f\|_{L^1}\|\tilde{\nabla} M^{(0)}_t\|_{L^1}.  
\end{equation*} 
Using the inequality (\ref{poly8}) we obtain
\begin{equation*} 
\int_{G}\left|t^{1-s/2}\mathcal{J}f\ast M^{(0)}_{t}(x)\right|dx\leq Ct^{\frac{1-s}{2}}\|\nabla f\|_{L^1}.
\end{equation*} 
\begin{flushright}{$\blacksquare$}\end{flushright}
\begin{Proposition}\label{PropoChame2}
For the last integral of (\ref{poids1}) we have the inequality
\begin{equation*} 
\int_{G}\left|t^{1-s/2}\mathcal{J}f\ast M^{(1)}_{t}(x)\right|dx\leq Ct^{\frac{1-s}{2}}\|\nabla f\|_{L^{1}} 
\end{equation*}
\end{Proposition} 
\emph{\textbf{Proof}}. Here, it is necessary to make an additional step. We cut out the function $m_{1}$ in the following way:  
\begin{equation*}
m_{1}(\lambda)=\left(\frac{1-e^{-\lambda}}{\lambda}\right)\theta_{1}(\lambda)=m_{a}(\lambda)-m_{b}(\lambda) 
\end{equation*} 
where $m_{a}(\lambda)=\frac{1}{\lambda}\theta_{1}(\lambda)$ and $m_{b}(\lambda)=\frac{e^{-\lambda}}{\lambda}\theta_{1}(\lambda)$. We will note $M^{(a)}_{t}$ and $M^{(b)}_{t}$ the associated kernels of these two operators. We obtain thus the estimate 
\begin{equation}\label{Rolling1} 
\int_{G}\left|t^{1-s/2}\mathcal{J}f\ast M^{(1)}_{t}(x)\right|dx\leq \int_{G}\left|t^{1-s/2}\mathcal{J}f\ast M^{(a)}_{t}(x)\right|dx+ 
\int_{G}\left|t^{1-s/2}\mathcal{J}f\ast M^{(b)}_{t}(x)\right|dx 
\end{equation} 
We have the next lemma for the last integral in (\ref{Rolling1}).
\begin{Lemme}\label{BishopLem0}
\begin{equation*} 
\int_{G}\left|t^{1-s/2}\mathcal{J}f\ast M^{(b)}_{t}(x)\right|dx\leq Ct^{\frac{1-s}{2}}\|\nabla f\|_{L^{1}}.
\end{equation*} 
\end{Lemme} 
\emph{\textbf{Proof}}. Observe that $m_{b}\in \mathcal{S}(\mathbb{R}^{+})$, then the proof is straightforward and follows the same steps as those of the preceding proposition \ref{Sir}.  
\begin{flushright}{$\blacksquare$}\end{flushright} 
We treat the other part of (\ref{Rolling1}) with the following lemma:  
\begin{Lemme} \label{BishopLem}
\begin{equation}\label{Bishop} 
\int_{G}\left|t^{1-s/2}\mathcal{J}f\ast M^{(a)}_{t}(x)\right|dx\leq Ct^{\frac{1-s}{2}}\|\nabla f\|_{L^{1}} 
\end{equation} \end{Lemme} 
\emph{\textbf{Proof}}. We consider the auxiliary function 
$$\psi(\lambda)=\theta_{0}(\lambda/2)-\theta_{0}(\lambda)=\theta_{1}(\lambda)-\theta_{1}(\lambda/2)$$ 
in order to obtain the identity $$\sum_{j=0}^{+\infty}\psi(2^{-j}\lambda)=\theta_{1}(\lambda).$$  
We have then 
$$m_{a}(t\lambda)=\frac{1}{t\lambda}\sum_{j=0}^{+\infty}\psi(2^{-j}t\lambda)=\sum_{j=0}^{+\infty}2^{-j}\phi(2^{-j}t\lambda)$$ 
where $\phi(\lambda)=\frac{\psi(\lambda)}{\lambda}$ is a function in $\mathcal{C}^{\infty}_{0}(\mathbb{R}^{+})$. Then, from the point of view of operators, one has:  
\begin{equation}\label{Stone} 
M^{(a)}_{t}=\sum^{+\infty}_{j=0}2^{-j}K_{j,t}
\end{equation} 
where $K_{j,t}=\phi(2^{-j}t\mathcal{J})$. With the formula (\ref{Stone}) we return to the left side of (\ref{Bishop}):  
\begin{equation}\label{Chicago}
\int_{G}\left|t^{1-s/2}\mathcal{J}f\ast M^{(a)}_{t}(x)\right|dx\leq\sum^{+\infty}_{j=0}2^{-j}\int_{G}\left|t^{1-s/2}\mathcal{J}f\ast 
K_{j,t}(x)\right|dx.
\end{equation} 
Using the sub-Laplacian definition and the vector fields properties, we have
$$\sum^{+\infty}_{j=0}2^{-j}\int_{G}\left|t^{1-s/2}\mathcal{J}f\ast K_{j,t}(x)\right|dx \leq \sum^{+\infty}_{j=0}2^{-j}t^{1-s/2}\|\nabla f\|_{L^1}\|\tilde{\nabla}K_{j,t}\|_{L^1}. $$
Apply now proposition \ref{toge} to obtain the estimate $\|\tilde{\nabla}K_{j,t}\|_{L^1}\leq C 2^{j/2}t^{-1/2}$. We have then for (\ref{Chicago}) the following inequality\\
$$\int_{G}\left|t^{1-s/2}\mathcal{J}f\ast M^{(a)}_{t}(x)\right|dx \leq C\sum^{+\infty}_{j=0}2^{-j/2}t^{\frac{1-s}{2}}\|\nabla f\|_{L^1}.$$ 
Then, we finally get 
\begin{equation*} 
\int_{G}\left|t^{1-s/2}\mathcal{J}f\ast M^{(a)}_{t}(x)\right|dx\leq C\, t^{\frac{1-s}{2} } \|\nabla f\|_{L^{1}}.  
\end{equation*} 
Which ends the proof of the lemma \ref{BishopLem}.
\begin{flushright}{$\blacksquare$}\end{flushright} 
With these two last lemmas we conclude the proof of the proposition \ref{PropoChame2}. Now, getting back to the formula (\ref{poids1}), with propositions \ref{Sir} and \ref{PropoChame2} we finally finish the proof of theorem \ref{CHAME}.
\begin{flushright}{$\blacksquare$}\end{flushright} 
\subsection{Weak inequalities}

To begin the proof, we use again the fact that operator $\mathcal{J}^{s/2}$ carries out an isomorphism between the spaces $\dot{B}^{-\beta, \infty}_{\infty}$ and $\dot{B}^{-\beta-s, \infty}_{\infty}$. Thus inequality (\ref{smile1}) rewrites as:  
\begin{equation}\label{ForIsoLap}
\|\mathcal{J}^{s/2}f\|_{L^{q, \infty}}\leq C \|\nabla f\|_{L^{1}}^{\theta}\|\mathcal{J}^{s/2}f\|_{\dot{B}^{-\beta-s, \infty}_{\infty}}^{1-\theta} 
\end{equation} 
By homogeneity, we can suppose that the norm $\|\mathcal{J}^{s/2}f \|_{\dot{B}^{-\beta-s,\infty}_{\infty}}$ is bounded by $1$; then we have to show 
\begin{equation}\label{CHAMEFAIBLE2} 
\|\mathcal{J}^{s/2}f \|_{L^{q, \infty}}\leq C \|\nabla f \|_{L^{1}}^{\theta}. 
\end{equation} 
We have thus to evaluate the expression $\left|\{x\in G:|\mathcal{J}^{s/2}f(x)|> 2\alpha\}\right|$ for all $\alpha>0$. If we use the thermic definition of the Besov space (\ref{BesovChaleurGhomo}), we have 
$$\|\mathcal{J}^{s/2}f\|_{\dot{B}^{-\beta-s,\infty}_{\infty}}\leq 1 \iff \underset{t>0}{\sup}\left\{t^{\frac{\beta+s}{2}}\|H_{t}\mathcal{J}^{s/2}f\|_{L^\infty}\right\} \leq 1.$$ 
But, if one fixes $t_{\alpha}=\alpha^{-\left(\frac{2}{\beta+s}\right)}$, we obtain $\|H_{t_{\alpha}}\mathcal{J}^{s/2}f \|_{L^\infty}\leq \alpha $. 
Note also that with the definition of parameter $\beta$ one has $t_{\alpha}=\alpha^{-\frac{2(q-1)}{(1-s)}}$. Therefore, since we have the following set inclusion
$$\left\{x\in G: |\mathcal{J}^{s/2}f(x)|> 2\alpha\right\}\subset \left\{x\in G: |\mathcal{J}^{s/2}f(x)-H_{t_{\alpha}}\mathcal{J}^{s/2}f(x)|> \alpha\right\}, $$ 
the Chebyshev inequality implies 
\begin{equation*} 
\alpha^{q}\left|\{x\in G: |\mathcal{J}^{s/2}f(x)|> 2\alpha\}\right|\leq \alpha^{q-1}\int_{G}|\mathcal{J}^{s/2}f(x)-H_{t_{\alpha}}\mathcal{J}^{s/2}f(x)|dx.  
\end{equation*} 
At this point, we use the theorem \ref{CHAME} to estimate the right side of the preceding inequality:  
\begin{equation}\label{T2g} 
\alpha^{q}\left|\{x\in G:|\mathcal{J}^{s/2}f(x)|> 2\alpha\}\right|\leq C \alpha^{q-1}\;t_{\alpha}^{\frac{1-s}{2}}\int_{G}|\nabla f(x)|dx.  
\end{equation} 
But, by the choice of $t_{\alpha}$, one has $\alpha^{q-1}\alpha^{-\frac{2(q-1)}{(1-s)}\frac{(1-s)}{2}}=1 $. Then (\ref{T2g}) implies the inequality 
$$\qquad\alpha^{q}\left|\{x\in G:|\mathcal{J}^{s/2}f(x)|> 2\alpha\}\right|\leq C \|\nabla f\|_ {L^{1}} \;; $$  
and, finally, using the definition (\ref{DefDebilSoboLe}) of weak Sobolev spaces it comes 
$$\qquad \qquad \|\mathcal{J}^{s/2}f \|^{q}_{L^{q, \infty}}\leq  C \|\nabla f \|_{L^{1}}$$ 
which is the desired result.\begin{flushright}{$\blacksquare$}\end{flushright} 
\subsection{Strong inequalities}
When $s=0$ in the weak inequalities above it is possible to obtain stronger estimations. To achieve this, we will need an intermediate step:  
\begin{Proposition}\label{escalera2} 
Let $1<q<+\infty$, $\theta=\frac{1}{q}$ and $\beta=\theta/(1-\theta)$.  Then we have
$$\|f\|_ { L^{q}}\leq C \|\nabla f\|_{L^{1}}^{\theta } \|f\|^{1-\theta}_{\dot{B}^{-\beta, \infty}_{\infty}}$$ 
when the three norms in this inequality are bounded.  
\end{Proposition} 
\textit{\textbf{Proof}.}  We will follow closely \cite{Ledoux}. Just as in the preceding theorem, we will start by supposing that $\|f\|_{\dot{B}^{-\beta, \infty}_{\infty}}\leq 1$.  Thus, we must show the estimate 
\begin{equation}\label{ForteG} 
\|f\|_{L^{q}}\leq C \|\nabla f\|_{L^{1}}^{\theta}. 
\end{equation} 
Let us fix $t$ in the following way:  $t_{\alpha}=\alpha^{-2(q-1)/q}$ where $\alpha>0$.  We have then, by the thermic definition of Besov spaces, the estimate $\|H_{t}f \|_{L^\infty}\leq \alpha$. We use now the characterization of Lebesgue space given by the distribution function:  
\begin{equation}\label{formulaG} 
\frac{1}{5^{q }} \|f\|_{L^{q}}^{q}=\int_{0}^{+\infty}\left|\{x\in G: |f(x)|> 5\alpha\}\right|d(\alpha^{q}).  
\end{equation} 
It now remains to estimate $|\{x\in G: |f(x)|> 5\alpha\}|$ and for this we introduce the following thresholding function:  
\begin{equation*}\label{teta}
\Theta_{\alpha}(t)=\left\lbrace\begin{array}{l} \Theta_{\alpha}(-t)=-\Theta_{\alpha}(t)\\[4mm] 
0 \qquad \qquad\qquad\mbox{ if }\qquad 0\leq T \leq \alpha\\[4mm]
t-\alpha \qquad\qquad\mbox{ if } \qquad\alpha\leq T \leq M\alpha \\[4mm]
(M-1)\alpha \qquad\mbox{ if } \qquad T > M\alpha 
\end{array}\right.\\[3mm ] 
\end{equation*}
Here, $M$ is a parameter which depends on $q$ and which we will suppose for the moment larger than 10.  \\

This cut-off function enables us to define a new function posing $f_{\alpha}=\Theta_{\alpha}(f)$. We collect in the next lemma some significant properties of this function $f_{\alpha}$: 
\begin{Lemme}\label{Taboo} 
\begin{enumerate} 
\item[]
\item The set defined by $\{x\in G: |f(x)|> 5\alpha\}$ is included in the set $\{x\in G: |f_{\alpha }(x)|> 4\alpha\}$. \\  
\item On the set $\{x\in G: |f(x)|\leq M\alpha\}$ one has the estimate $|f-f_{\alpha }|\leq \alpha$. \\  
\item If $f\in \mathcal{C}^{1}(G)$, one has the equality $\nabla f_{\alpha}=(\nabla f)\mathds{1}_{\{\alpha\leq|f|\leq M\alpha\}}$ 
almost everywhere.  
\end{enumerate} 
\end{Lemme} 
We leave the verification of this lemma to the reader. \\  

Let us return now to (\ref{formulaG}). By the first point of the lemma above we have
\begin{equation}\label{cray} 
\int_{0}^{+\infty}\left|\left\{x\in G: |f(x)|> 5\alpha\right\}\right|d(\alpha^{q})\leq \int_{0}^{+\infty }\left|\{x\in G: |f_{\alpha }(x)|> 4\alpha\}\right|d(\alpha^{q})=I. 
\end{equation} 
We note $A_{\alpha}=\{x\in G: |f_{\alpha }(x)|> 4\alpha\}$, $B_{\alpha}=\{x\in G: |f_{\alpha}(x)-H_{t_{\alpha}}(f_{\alpha})(x)|> \alpha\}$ and $C_{\alpha}=\{x\in G: |H_{t_{\alpha}}(f_{\alpha}-f)(x)|> 2\alpha\}$. Now, by linearity of $H_{t}$ we can write: $f_{\alpha}=f_{\alpha}-h_{t_{\alpha}}(f_{\alpha})+h_{t_{\alpha}}(f_{\alpha}-f)+h_{t_{\alpha}}(f)$. Then, holding in account the fact $\|H_{t}f\|_{L^\infty}\leq \alpha$, we obtain $A_{\alpha}\subset B_{\alpha}\cup C_{\alpha}$. Returning to (\ref{cray}), this set inclusion gives us the following inequality 
\begin{equation}\label{dosG} 
I\leq \int_{0}^{+\infty}\left|B_{\alpha}\right|d(\alpha^{q})+\int_{0}^{+\infty}\left|C_{ \alpha}\right|d(\alpha^{q }) 
\end{equation} 
We will study and estimate these two integrals, which we will call $I_{1}$ and $I_{2}$ respectively, by the two following lemmas:  
\begin{Lemme} For the first integral of (\ref{dosG}) we have the estimate:  
\begin{equation}\label{ARG} 
I_{1 } = \int_{0}^{+\infty}\left|B_{\alpha}\right|d(\alpha^{q})\leq C\, q\log(M) \|\nabla f\|_ { L^{1} } 
\end{equation} 
\end{Lemme} 
\textbf{\emph{Proof}.} The Chebyshev's inequality implies 
\begin{equation*} 
\left|B_{\alpha}\right|\leq \alpha^{-1}\int_{G}|f_{\alpha}(x)-H_{t_{\alpha}}(f_{\alpha })(x)|dx.  
\end{equation*} 
Using the theorem \ref{CHAME} with $s=0$ in the above integral we obtain:  
$$\left|B_{\alpha}\right|\leq C \, \alpha^{-1 } \, t_{\alpha}^{1/2}\int_{G}|\nabla f_{\alpha }(x)|dx$$ 
Remark that the choice of $t_{\alpha}$ fixed before gives $t_{\alpha}^{1/2}=\alpha^{1-q}$, then we have
$$\left|B_{\alpha}\right|\leq C \, \alpha^{-q}\int_{\{\alpha\leq|f|\leq M\alpha \} }|\nabla f(x)|dx.$$ 
We integrate now the preceding expression with respect to $d(\alpha^{q})$:  
$$I_{1}\leq C\int_{0}^{+\infty}\alpha^{-q}\left(\int_{\{\alpha\leq|f|\leq M\alpha \} }|\nabla f(x)|dx\right)d(\alpha^{q }) = C\;q\int_{G}|\nabla f(x)|\left(\int_{\frac{|f|}{M}}^{|f|}\frac{d\alpha}{\alpha}\right)dx$$ 
It follows then $I_{1}\leq C\, q\, \log(M) \|\nabla f \|_{L^{1}}$ and one obtains the estimation needed for the first integral.
\begin{flushright}{$\blacksquare$}\end{flushright} 
\begin{Lemme}\label{papas}
For the second integral of (\ref{dosG}) one has the following result:  
\begin{equation*}
I_{2}=\int_{0}^{+\infty}\left|C_{\alpha}\right|d(\alpha^{q})\leq \frac{q}{q-1}\;\frac{1}{M^{q-1}} \|f\|_{L^q}^{q} 
\end{equation*} 
\end{Lemme} 
\textbf{\emph{Proof}.} For the proof of this lemma, we write:  $$|f-f_{\alpha }|=|f-f_{\alpha }|\mathds{1}_{\{|f|\leq M \alpha\}}+|f-f_{\alpha }|\mathds{1}_{\{|f|> M\alpha\}}.$$ 
As the distance between $f$ and $f_{\alpha}$ is lower than $\alpha$ on the set $ \{x\in G: |f(x)|\leq M \alpha\}$, one has the inequality 
$$|f-f_{\alpha }|\leq\alpha+|f|\mathds{1}_{\{|f|> M \alpha\}}$$ 
By applying the heat semi-group to both sides of this inequality we obtain $H_{t_{\alpha}}(|f-f_{\alpha }|)\leq \alpha + H_{t_{\alpha}}(|f|\mathds{1}_{\{|f|> M \alpha\}})$ and we have then the following set inclusion $C_{\alpha}\subset \left\{x\in G: H_{t_{\alpha}}(|f|\mathds{1}_{\{|f|> M \alpha\}})>\alpha\right\}$.
Thus, considering the measure of these sets and integrating with respect to $d(\alpha^{q})$, it comes 
\begin{equation*} 
I_{2}=\int_{0}^{+\infty}\left|C_{\alpha}\right|d(\alpha^{q})\leq\int_{0}^{+\infty}\bigg|\{H_{t_{\alpha}}(|f|\mathds{1}_{\{|f|> M\alpha\}})>\alpha\}\bigg|d(\alpha^{q}) 
\end{equation*} 
We obtain now, by applying the Chebyshev inequality, the estimate 
$$I_{2}\leq\int_{0}^{+\infty}\alpha^{-1}\bigg(\int_{G}H_{t_{\alpha}}\left(|f|\mathds{1}_{\{|f|> M\alpha\}}\right)dx\bigg)d(\alpha^{q}),$$ 
then by Fubini's theorem we have
$$I_{2}\leq q\int_{G}|f(x)|\bigg(\int_{0}^{+\infty}\mathds{1}_{\{|f|> M\alpha\}}\alpha^{q-2}d\alpha\bigg)dx=\frac{q}{q-1}\int_{G}|f(x)|\frac{|f(x)|^{q-1}}{M^{q-1}}dx= \frac{q}{q-1}\frac{1}{M^{q-1 }}\|f\|_{L^{q}}^{q}.$$ 
And this concludes the proof of this lemma.
\begin{flushright}{$\blacksquare$}\end{flushright} 
We finish the proof of proposition \ref{escalera2} by connecting together these two lemmas \textit{i.e.}:
$$\frac{1}{5^{q}} \|f\|_ {L^q}^{q}\leq Cq\, \log(M) \|\nabla f\|_ {L^1}+\frac{q}{q-1}\frac{1}{M^{q-1 } } \|f\|_{L^q}^{q}$$ 
Since we supposed all the norms bounded and  $M\gg 1$, we finally have
$$\left(\frac{1}{5^{q}}-\frac{q}{q-1}\frac{1}{M^{q-1}}\right) \|f\|_ {L^q}^{q}\leq C q\, \log(M) \|\nabla f\|_{L^1}$$ 
\begin{flushright}{$\blacksquare$}\end{flushright} 
The proof of the theorem \ref{smile00} is not yet completely finished. The last step is provided by the 
\begin{Proposition}\label{escalera3} 
In the proposition \ref{escalera2} it is possible to consider only the two assumptions $\nabla f\in L^{1}(G)$ and $f\in \dot{B}^{-\beta, \infty}_{\infty}(G)$.  
\end{Proposition}
\textbf{\emph{Proof}.} For the proof of this proposition we will build an approximation of $f$ writing:  
\begin{equation*} 
f_{j}=\left(\int_{0}^{+\infty}\left(\varphi(2^{-2j}\lambda)-\varphi(2^{2j}\lambda)\right)dE_{\lambda}\right)(f)
\end{equation*} 
where $\varphi$ is a $\mathcal{C}^{\infty}(\mathbb{R}^+)$ function such that $\varphi=1$ on $]0,1/4[$ and $\varphi=0$ on $[1, +\infty[$.
\begin{Lemme} 
If $q>1$, if $\nabla f\in L^{1}(G)$ and if $f\in \dot{B}^{-\beta, \infty}_{\infty}(G)$ then $\nabla f_j\in L^{1}(G)$, $f_j\in \dot{B}^{-\beta, \infty}_{\infty}(G)$ and $f_{j}\in L^{q}(G)$.  
\end{Lemme} 
\textbf{\emph{Proof}.} The fact that $\nabla f_j\in L^{1}(G)$ and $f_j\in \dot{B}^{-\beta, \infty}_{\infty}(G)$ is an easy consequence of the definition of $f_j$. For $f_{j}\in L^{q}(G)$ the starting point is given by the relation:  
\begin{equation*}
f_{j}=\left(\int_{0}^{+\infty}m(2^{-2j}\lambda)\, dE_{\lambda}\right)2^{-2j}\mathcal{J}(f), 
\end{equation*} 
where we noted $$m(2^{-2j}\lambda)=\frac{\varphi(2^{-2j}\lambda)-\varphi(2^{2j}\lambda)}{2^{-2j}\lambda}.$$ 
Observe that the function $m$ vanishes near of the origin and satisfies the assumptions of proposition \ref{toge}. We obtain then the following identity where $M_{j}$ is the kernel of the operator $m(2^{-2j}\mathcal{J})$:  
$$f_{j}=2^{-2j}\mathcal{J}f\ast M_{j}=2^{-2j}\nabla f\ast \tilde{\nabla}M_{j},$$ 
Using inequality (\ref{poly8}), we estimate the norm $L^{q}(G)$ in the preceding identity:
$$\|f_{j}\|_{L^{q}}=\|2^{-2j}\nabla f\ast\tilde{\nabla}M_{j}\|_{L^{q}}\leq 2^{-2j}\|\nabla f\|_{L^{1}}\|\tilde{\nabla}M_{j}\|_{L^{q}}.$$ 
Finally, we obtain:  
$$\|f_{j} \|_{L^{q}}\leq C\, 2^{j(d(1-\frac{1}{q})-1) } \|\nabla f \|_{L^{1}}<+\infty$$ 
\begin{flushright}{$\blacksquare$}\end{flushright} 

Thanks to this estimate, we can apply the proposition \ref{escalera2} to $f_{j}$ whose $L^{q}(G)$ norm is bounded, and we obtain:  
$$\|f_{j}\|_{L^{q}}\leq C \|\nabla f_j\|_{L^{1}}^{\theta } \|f_j\|_{\dot{B}^{-\beta, \infty}_{\infty}}^{1-\theta}.$$ 
Now, since $f\in \dot{B}^{-\beta, \infty}_{\infty}(G)$, we have $f_{j} \rightharpoonup f$ in the sense of distributions. It follows  
$$\|f\|_{ L^{q}}\leq \underset{j \to +\infty}{\lim \inf } \|f_{j } \|_ { L^{q}}\leq C \|\nabla f\|_ { L^{1}}^{\theta }
 \|f\|_ { \dot{B}^{-\beta, \infty}_{\infty}}^{1-\theta}.$$ 
We restricted ourselves to the two initial assumptions, namely $\nabla f\in L^{1}(G)$ and $f\in\dot{B}^{-\beta, \infty}_{\infty}(G)$. The strong inequalities (\ref{smile0}) are now completely proved for stratified groups.  
\begin{flushright}{$\blacksquare$}\end{flushright} 


\begin{flushright}
\begin{minipage}[r]{80mm}
Diego \textsc{Chamorro}\\[5mm]
Laboratoire d'Analyse et de Probabilités\\ 
Université d'Evry Val d'Essonne \& ENSIIE\\[2mm]
1 square de la résistance,\\
91025 Evry Cedex\\[2mm]
diego.chamorro@m4x.org
\end{minipage}
\end{flushright}

\end{document}